\documentclass[11pt,epsfig]{article}
\usepackage{amsmath}
\usepackage{amsthm}
\usepackage{epsfig}
\usepackage{graphicx}
\usepackage{hyperref}
\usepackage{amsfonts}

\newtheoremstyle{thmm}{1.5ex plus 1ex minus .2ex}{1.5ex plus 1ex minus .2ex}{\rmfamily}{}{\bfseries}{}{1em}{}
\theoremstyle{thmm}
\newtheorem{theorem}{Theorem}[section]
\newtheorem{lemma}{Lemma}[section]

\newtheorem{definition}{Definition}[section]

\renewenvironment{proof}[1][Proof]{\noindent\textit{#1. } }{\hfill$\square$}

\newcommand{\nn}{\nonumber}
\def \endproof{\vrule height8pt width 5pt depth 0pt}

\def\v{\varepsilon}

\def\t{\theta}

\def\k{\kappa}

\def\a{\alpha}
\def\b{\beta}

\def\G{\Gamma}

\def\l{\lambda}
\def\f{\frac}

\def\r{\rho}

\def\s{\sigma}

\def\di{\displaystyle}
\def\i{\infty}
\def\R{\mathbb{R}}

\date{}
\title{\bf  Global existence of weak solution
to the heat and moisture transport system in fibrous porous media}
\author{Buyang Li \footnote{Department of Mathematics, City University of Hong Kong,
Kowloon, Hong Kong ({\tt buyangli2@student.cityu.edu.hk,
maweiw@math.cityu.edu.hk}). The work of the authors was supported
in part by a grant from the Research Grants Council of the Hong
Kong Special Administrative Region, China (Project No. CityU
102005).}  \ , Weiwei Sun \footnotemark[1] \ , and Yi Wang
\footnote{ Institute of Applied Mathematics, AMSS, CAS, Beijing
100190, China ({\tt wangyi@amss.ac.cn}). The work of this author
was supported in part by the NSFC grant (No. 10801128) and a grant
from the Research Grants Council of the Hong Kong Special
Administrative Region, China (Project No. CityU 102005). } }

\begin{document}
\maketitle

\begin{abstract}
This paper is concerned with theoretical analysis of
a heat and moisture transfer model arising from textile
industries, which is described by a degenerate and strongly
coupled parabolic system. We prove the global (in time)
existence of weak solution by constructing
an approximate solution with some standard smoothing.
The proof is based on the physcial nature of gas convection, in which
the heat (energy) flux in convection is determined by
the mass flux in convection.

\end{abstract}


\section{Introduction}
\setcounter{equation}{0} Mathemaitcal modeling for heat and moisture
transport with phase change in porous textile materials was studied
by many authors, $e.g.$, see \cite{FLL,FCWS,HYS,LZ}. A typical
application of these models is a clothing assembly, consisting of a
thick porous fibrous batting sandwiched by two thin fabrics. The
outside cover of the assembly is exposed to a cold environment with
fixed temperature and relative humidity while the inside cover is
exposed to a mixture of air and vapor at higher temperature and
relative humidity. In general, the physical process can be viewed as
a multiphase and single (or multi) component flow. In this process,
the water vapor moves through the clothing assembly by convection
which is induced by the pressure gradient. The heat is transferred
by conduction in all phases (liquid, fiber and gas) and convection
in gas. Phase changes occur in the form of evaporation/condensation
and/or sublimation. Based on the conservation of mass and energy and
the neglect of the water influence, the model can be described by
\begin{eqnarray}
& &\frac{\partial}{\partial t}(\epsilon
C_v)+\frac{\partial}{\partial x}(u \epsilon C_v) =-\Gamma_{ce},
\label{Cv-e}\\
&&\frac{\partial }{\partial t} \left ( \epsilon C_{vg}M C_vT +
(1-\epsilon) C_{vs} T \right ) + \frac{\partial}{\partial x} \left
( \epsilon C_{vg} M u C_vT \right ) = \frac{\partial}{\partial x}
\left ( \kappa \frac {\partial T}{\partial x} \right ) + \lambda
M\Gamma_{ce} \, . \label{T-e}
\end{eqnarray}
Here $C_v$ is the vapor concentration ($\rm mol/m^3$), $T$ is the
temperature ($K$), $\epsilon$ the porosity of the fiber,
$M$ the molecular weight of water and $\lambda$  the latent heat of
evaporation/condensation in the wet zone while in frozen zone, it
represents the latent heat of sublimation. $C_{vg}$ and $C_{vs}$
are the heat capacities of the gas and mixture solid,
respectively.

The evaporation/condensation (molar) rate of phase change per unit
volume is defined by the Hertz-Knudsen equation~\cite{Jon}
\begin{equation}
\Gamma_{ce} = -\frac{E}{R_f}
\sqrt{\frac{(1-\epsilon)(1-\epsilon')}{2\pi RM }} \left(\frac{P_{\rm
sat}(T)- P}{\sqrt{T}}\right) \label{phase}
\end{equation}
where $R$ is the universal gas constant, $R_f$ is the radius of
fibre and $E$ is the nondimensional phase change coefficient. The
vapor pressure is given by $P=RC_vT$ because of the ideal gases'
assumption. The saturation pressure $P_{\rm sat}$ is determined from
experimental measurements, see Figure \ref{sat}.


The vapor velocity (volumetric discharge) is given by the Darcy's
law
\begin{equation}
u=-\frac{k k_{rg}}{\mu_{g}}\frac{\partial P}{\partial x
}\label{darcy-v}
\end{equation}
where $k$ is the permeability, $k_{rg}$ and $\mu_{g}$ are the
relative permeability and the viscosity of the vapor,
respectively.

Numerical methods and simulations for the heat and moisture
transport in porous textile materials have been studied by many
authors with various applications \cite{CW,HSY,OT,ST,WMR}. However,
no theoretical analysis has been explored for the above system of
nonlinear equations. A simple heat and moisture model was studied
in \cite{Val}, where the model was described by a pure diffusion
process (without convection and condensation) with a non-symmetric
parabolic part. There are several related porous media flow
problems from other physical applications. A popular one is a
compressible (or incompressible) flow in porous
media with applications in oil and underground water industries,
which is described by an elliptic pressure equation coupled with a
parabolic concentration equation for incompressible case and a
system of parabolic equations for compressible case. The existence
of weak solution for incompressible and compressible flows has been
studied in \cite{CE, Feng} and \cite{AS}, respectively. However, in
most of these works, the temperature is ignored and the phase change
(condensation/evaporation) does not occur due to the nature of these
applications, while both temperature and phase change play important
roles in the textile model.

For the textile model, the water content in the batting area
usually is relative small and one often assumes that all these
physical parameters involved in the system
(\ref{Cv-e})-(\ref{T-e}) are positive constants. With
nondimensionalization, the system (\ref{Cv-e})-(\ref{T-e}) reduces
to
\begin{equation}
\left\{
\begin{array}{l}
\r_t-((\rho \theta)_x \r)_x=-\G,  \\
(\rho \theta)_t+\s\t_t-((\rho \theta)_x \rho \theta
)_x-(\k\t_x)_x=\l\G, ,
\end{array}
\right.\label{sys}
\end{equation}
for $x \in (0,1)$, $t>0$, where $(\cdot)_{\mu} =
\frac{\partial}{\partial \mu}$ for $\mu=x, t$, $\r=\r(x,t)$ and
$\t=\t(x,t)$ represent the density of vapor and the temperature,
respectively,
$$
\Gamma = \rho \theta^{1/2} - p_s(\theta) \,
$$
and $p_s(\theta) \sim P_{\rm sat}(\theta)/\theta^{1/2}$. $\s$ and
$\l$ are given positive constants and $\k=\k_1+\k_2\r^2$ is the heat
conductivity coefficient with $\k_i~(i=1,2)$ being positive
constants. We consider a class of commonly used Robin type boundary
conditions \cite{FLL,FCWS,HYS,YHFS} defined by
\begin{equation}
(\rho \theta)_x\r|_{x=1}=\a^1(\bar\r^1-\r(1,t)),\quad (\rho
\theta)_x\r|_{x=0}=\a^0(\r(0,t)-\bar\r^0), \label{sys-b1}
\end{equation}
and
\begin{equation} \k\t_x|_{x=1}=\b^1(\bar\t^1-\t(1,t)),\quad
\k\t_x|_{x=0}=\b^0(\t(0,t)-\bar\t^0),\label{sys-b2}
\end{equation}
and the initial condition is
\begin{equation}
\r(x,0)=\r_0(x),\quad\t(x,0)=\t_0(x),\qquad x\in (0,1)
\label{sys-i}
\end{equation}
where $\a^0,\a^1$ represent the mass transfer coefficients,
$\bar\r^0,\bar\r^1$ are the density of the gas in the inner
background and outer background, respectively, $\b^0,\b^1$ the
heat transfer coefficients, and $\bar\t^0,\bar\t^1$ the inner and
outer background temperatures. Physically, all the parameters
above are positive constants and
$\r_0(x)>\underline{\r},\t_0(x)>\underline{\t}$, where
$\underline{\r}$ and $\underline{\t}$ are positive constants.
Based on the experimental data in Figure \ref{sat}, we assume that
$p_s$ is a smooth, increasing and nonnegative function defined on
$\R^+$ which satisfies
\begin{equation}
\lim_{\theta \rightarrow 0 }
\frac{p_s(\theta)}{\theta}=0,\;\;\;\quad \lim_{\theta \rightarrow
\infty } \frac{p_s(\theta)}{\theta^{1+\eta}}=\infty \label{sat-a}
\end{equation}
for some $\eta>0$. For physical reasons, we set $p_s(\t)=0$ for
$\t\leq 0$.

The objective of this paper is to establish the global existence of
weak solution to the initial-boundary value problem
(\ref{sys})-(\ref{sys-i}) under the general physical hypotheses
(\ref{sat-a}) for the saturation pressure function $\Gamma$. The
difficulty lies on the strong nonlinearity and the coupling of
equations. To the best of our knowledge, there are no theoretical
results for the underlying model. More important is its significant
applications in textile industries. Also analysis presented in this
paper may provide a fundamental tool for theoretical analysis of
existing numerical methods. Our proof is based on the equivalence of
mass and heat transfer in convection.

\section{The main result}
\setcounter{equation}{0} Before we present our main result, we
introduce some notations. Let $T$ be a given positive number in
the following sections. We define
$$
\Omega=(0,1),\quad I=(0,T],\quad Q_t = \Omega\times (0,t], \quad
Q_T= \Omega\times I,
$$
$$
V_1(Q_T)=L^2(I;H^1(\Omega)), \quad V_2(Q_T)=\Big\{f\in
L^2(Q_T)~\Big|~\|f\|_{V_2(Q_T)}<+\i\Big\},
$$
$$
\|f\|_{V_2(Q_T)}= {\rm ess}\sup_{\!\!\!\!t\in [0,T]}
\|f(\cdot,t)\|_{L^2(\Omega)}+\|f_x\|_{L^2(Q_T)},
$$
$$
W_2^{2,1}(Q_T)=\Big\{f\in L^2(Q_T)~\Big|~f_t,f_x,f_{xx} \in
L^2(Q_T)\Big\}.
$$
Let $\mathcal{D}(\overline\Omega\times[0,T))$ be the subspace of
$C^\infty(\R^2)$ consisting of functions which have compact
support in $\R\times[-\infty,T)$, restricted to
$\overline\Omega\times[0,T)$.

Now we give the definition of weak solution to the system
(\ref{sys})-(\ref{sys-i}) and then, state our main result.
\vskip0.1in

\begin{definition}\label{defweaksol}$\!\!$ {\bf(Weak solution)}$\;$
{\it We say that the measurable function pair $(\r,\t)$ defined on
$\overline\Omega\times[0,T)$ is a global weak solution to
(\ref{sys})-(\ref{sys-i}) if $(\r,\t)\in (V_1(Q_T))^2$ and the
density $\r$ and the temperature $\t$ are nonnegative functions
satisfying
\begin{align}\label{rdefeq}
&\int_0^T\alpha^0(\r(0,t)-\bar\r^0)\phi(0,t)dt
+\int_0^T\alpha^1(\r(1,t)-\bar\r^1)\phi(1,t)dt
\nn\\
&+\int_0^T\int_\Omega (-\r\phi_t+(\r\t)_x\r \phi_x+\G\phi) dxdt
=\int_\Omega\rho_0\phi_0dx
\end{align}
and
\begin{align}\label{tdefeq}
&\int_0^T[\alpha^0(\r(0,t)-\bar\r^0)\t(0,t)+\beta^0(\t(0,t)-\bar\t^0)]
\psi(0,t)dt
\nn\\
&~~+\int_0^T[\alpha^1(\r(1,t)-\bar\r^1)\t(1,t)
+\beta^1(\t(1,t)-\bar\t^1)]\psi(1,t)dt
\nn\\
&~~+\int_0^\infty\int_\Omega \Big[-(\r\t+\s\t)\psi_t+(\r\t)_x\r\t
\psi_x+\k\t_x\psi_x-\l\G\psi\Big] dxdt
\nn\\
&=\int_\Omega (\r_0\t_0+\s\t_0)\psi_0dx
\end{align}
for any test functions $\phi,\psi\in
\mathcal{D}(\overline\Omega\times[0,T))$. }
\end{definition}
\vskip0.1in

\begin{theorem}\label{thm1}
{\it If the initial value $(\r_0,\t_0)$ satisfies $\r_0\in
L^{1+\gamma}(\Omega)$ $(\gamma>0)$, $\t_0\in L^\i(\Omega)$ and
$\rho_0 \ge 0$, $\t_0\geq\underline{\t}$ for some positive
constant $\underline{\t}$, then there exists a global weak
solution $(\r,\t)$, in the sense of Definition \ref{defweaksol},
to the initial-boundary value problem (\ref{sys})-(\ref{sys-i})
such that
$$
\r\ln \r\in L^\i(0,T;L^1(\Omega)),\quad \r\in
L^4(Q_T),\quad\r_x\in L^2(Q_T);
$$
$$
\t,\t^{-1}\in L^\i(Q_T),\quad  (1+\r)\t_x\in L^2(Q_T).
$$
}
\end{theorem}

In the following sections, we denote by $C_{p_1,p_2,\cdots,p_k}$ a
generic positive constant, which depends solely upon
$p_1,p_2,\cdots,p_k$, the physical parameters $\kappa_1, \kappa_2,
\sigma$ and $\lambda$ and the parameters involved in initial and
boundary conditions. In addition, we denote by
$C(p_1,p_2,\cdots,p_k)$ a generic positive function, dependent
upon the physical parameters $\kappa_1, \kappa_2, \sigma$ and
$\lambda$ and the parameters involved in boundary conditions,
which is bounded when $p_1,p_2,\cdots,p_k$ are bounded.

\section{Construction of approximate solutions}\label{apprxsol}
\setcounter{equation}{0} Throughout this section, we let $\v$ be a
fixed positive number which satisfies
$$
0<\v\leq\min\{\bar\r^0,\bar\r^1,\bar\t^0,\bar\t^1,1\},
$$
and $0<\nu<\v$. To prove the existence of global weak solutions to
the system (\ref{sys})-(\ref{sys-i}), we introduce a regularized
approximate system as follows:
\begin{align}
&\r_t-((\v+(\r\t)_\nu)\r_x)_x-(\r(\r_\v\t_x)_\v)_x
=-\r\chi^\v(\sqrt{\t})+\chi^\v(p_s(\t)),
\nn \\[3mm]
&(\r\t+\s\t)_{t}-(\kappa^\v\t_x)_x-((\v+(\r\t)_\nu))
\r_x\t)_x-(\r(\r_\v\t_x)_\v\t)_x \label{asys2}
\\
&=\l \r\chi^\v(\sqrt{\t})-\l\chi^\v(p_s(\t)) +(\l+\t) \left (
\chi^\v(p_s(\t))-p_s(\t)\right ) , ~~{\rm in}~~~ Q_T,\nn
\end{align}
where $\chi^\v$ is a cut-off function defined by
\begin{equation}
\chi^\v(h)=\left\{
\begin{array}{lll}
h&\mbox{\rm if}& |h|\leq \v^{-1},\\[3mm]
\v^{-1}&\mbox{\rm if}& |h|\geq \v^{-1},
\end{array}
\right. \nn
\end{equation}
and
$$
\kappa^\v=\kappa_1+\kappa_2|\r_\v|^2,
$$
and the subscriptions $\v, \nu$ define the smoothing operators in
general by $f_{\mu} = {\rm Ext}(f)*\eta_{\mu}$ with $\mu=\nu, \v$.
Here $\eta_{\mu}$ is the standard mollifier and ${\rm Ext}(\cdot)$
is the extension operator which extends any measurable functions
defined on $\Omega_T$ to be zero on $\R^2\backslash\Omega_T$.

The system (\ref{asys2}) can be rewritten as
\begin{equation}
\left\{
\begin{array}{l}
\r_t-((\v+(\r\t)_\nu)\r_x)_x-(\r(\r_\v\t_x)_\v)_x
+\r\chi^\v(\sqrt{\t})=\chi^\v(p_s(\t)),\\[5pt]
(\r+\s)\t_t-(\k^\v\t_x)_x-\left[(\v+(\r\t)_\nu)\r_x+\r(\r_\v\t_x)_\v\right]
\t_x -\r\chi^\v(\sqrt{\t})\t+(\l+\t)
p_s(\t)\di=\l\r\chi^\v(\sqrt{\t}) \, .
\end{array}
\right. \label{asys}
\end{equation}
The corresponding initial and  boundary conditions are given by
\begin{equation}
\begin{array}{l}
\di (\v+(\r\t)_\nu)\r_x+\r(\r_\v\t_x)_\v\big|_{x=1}
=\a^1(\bar\r^1-\r(1,t)),\\[3mm]
\di (\v+(\r\t)_\nu)\r_x+\r(\r_\v\t_x)_\v\big|_{x=0}
=\a^0(\r(0,t)-\bar\r^0),\\[3mm]
\di \r(x,0)=\r_{0\v}(x):=(\r_0)_\v(x)+\v,\\[3mm]
\di \k^\v\t_x|_{x=1}=\b^1(\bar\t^1-\t(1,t)),\\[3mm]
\di \k^\v\t_x|_{x=0}=\b^0(\t(0,t)-\bar\t^0),\\[3mm]
\di \t(x,0)=\t_{0\v}(x):=(\t_0)_\v(x).
\end{array}
\label{asys-bi}
\end{equation}

We prove the existence of solutions to the system
(\ref{asys})-(\ref{asys-bi}) by using the Leray--Schauder fixed
point theorem. The following lemma (see \cite{Lions1},
\cite{Lions}) is useful in our proof. \vskip0.1in

\begin{lemma}
$\!\!${\bf(Aubin--Lions)}$\;$ {\it Let
$B_1\hookrightarrow\hookrightarrow B_2\hookrightarrow B_3$ be
reflective and separable Banach spaces. Then
$$\{u\in L^p(I;B_1)|\;u_t\in
L^q(I;B_3)\}\hookrightarrow\hookrightarrow L^p(I;B_2), \quad
1<p,q<\infty;$$
$$\{u\in L^q(I;B_2)\cap L^1(I;B_1)|\;u_t\in L^1(I;B_3)\}
\hookrightarrow\hookrightarrow L^p(I;B_2),\quad 1\leq
p<q<\infty.$$ }
\end{lemma}

\subsection{Existence of approximate solutions}
We define
\begin{align*}
X=\{u\in L^2(I;H^1(\Omega))|\;u\geq0\},\quad Y=\{u\in
W^{2,1}_2(Q_T)|\;u\geq0\}.
\end{align*}
By Aubin--Lions lemma, $Y\hookrightarrow\hookrightarrow X$. Let $\v$
and $\nu$ be given positive constants and the parameter $s\in[0,1]$.
For any given $(\r^0,\t^0)\in X^2$, we define $\r$ to be the
solution of the following linear parabolic equation
\begin{equation}
 \r_t-((\v+(\r^0\t^0)_\nu)\r_x)_x-(\r(\r^0_\v\t^0_x)_\v)_x
+s\r\chi^\v(\sqrt{\t^0})=s\chi^\v(p_s(\t^0)), \label{rho-e}
\end{equation}
with the initial and boundary conditions
\begin{equation}
\left\{
\begin{array}{lr}
\di(\v+(\r^0\t^0)_\nu)\r_x+\r(\r^0_\v\t^0_x)_\v
=\a^1(s\bar\r^1-\r),&\mbox{\rm at}~~~ x=1,
\\[3mm]
\di(\v+(\r^0\t^0)_\nu)\r_x+\r(\r^0_\v\t^0_x)_\v
=\a^0(\r-s\bar\r^0),&\mbox{\rm at}~~~
 x=0,\\[3mm]
\r(x,0)=s\r_{0\v}(x),&\mbox{\rm for}~~~x\in\Omega.
\end{array}
\right. \label{rho-b}
\end{equation}
Now with $\r$ in hand, we define $\t$ to be the solution of the
linear parabolic equation
\begin{equation}
\begin{array}{r}
(\r+\s)\t_t-(\k^\v\t_x)_x-\big[(\v+(\r^0\t^0)_\nu)\r_x
+\r(\r^0_\v\t^0_x)_\v\big]\t_x\\[6pt]
-s\r\chi^\v(\sqrt{\t^0})\t+s(\lambda+\t)
p_s(\t)\di=s\lambda\r\chi^\v(\sqrt{\t^0}),
\end{array}
\label{theta-e}
\end{equation}
with the initial and boundary conditions
\begin{equation}
\left\{
\begin{array}{lr}
\di \k^\v\t_x=\b^1(s\bar\t^1-\t),&\mbox{\rm at}~~~
 x=1,\\[3mm]
\di \k^\v\t_x=\b^0(\t-s\bar\t^0),&\mbox{\rm at}~~~
 x=0,\\[3mm]
\t(x,0)=s\t_{0\v}(x),&\mbox{for}~~~x\in\Omega \, .
\end{array}
\right.\label{theta-b}
\end{equation}
Let $M$ denote the mapping from $(\r^0,\t^0,s)$ to $(\r,\t)$. Then
we have the following lemma. \vspace{5pt}

\begin{lemma}\label{mpcompconti}
{\it The mapping $M:X^2\times[0,1]\rightarrow X^2$ is well
defined, continuous and compact. }
\end{lemma}

{\it Proof}. By the $L^2$-theory of linear parabolic equations
\cite{LSU}, there exists a solution $\r\in W^{2,1}_2(Q_T)$ for the
system (\ref{rho-e})-(\ref{rho-b}) such that
$$
\|\r\|_{W^{2,1}_2(Q_T)}\leq
C(\v^{-1},\|(\r^0\t^0)_\nu\|_{C^1(\overline
Q_T)},\|(\r^0_\v\t^0_x)_\v)_x\|_{C^1(\overline
Q_T)},\|\r_{0\v}\|_{H^1(\Omega)},T).
$$
By noting the fact
$$
\|\r^0_\v\|_{H^1(\Omega)}\leq C_\v\|\r_0\|_{L^1(\Omega)}, \quad
\|(\r^0\t^0)_\nu\|_{C^1(\overline Q_T)}\leq
C_{\nu,T}\|\r^0\|_{L^2(Q_T)}\|\t^0\|_{L^2(Q_T)},
$$
$$\|(\r^0_\v\t^0_x)_\v)_x\|_{C^1(\overline
Q_T)}\leq C_{\v,T}\|\r^0_\v\t^0_x\|_{L^1(Q_T)}\leq
C_{\v,T}\|\r^0\|_{L^2(Q_T)}\|\t^0_x\|_{L^2(Q_T)} \, ,
$$
for the standard smoothing operator, we have
\begin{align}\label{rhomapping1}
\|\r\|_{W^{2,1}_2(Q_T)}\leq
C(\v^{-1},\nu^{-1},\|\r^0\|_X,\|\t^0\|_X,T)
\end{align}
and therefore,
$$
\|\r\|_{L^\infty(Q_T)}\leq \|\r\|_{W^{2,1}_2(Q_T)}\leq
C(\v^{-1},\nu^{-1},\|\r^0\|_X,\|\t^0\|_X,T).
$$

Let $\r^+=\max\{\r,0\}$, $\r^-=\max\{-\r,0\}$. Then $\r=\r^+-\r^-$.
By multiplying $\r^-$ on both sides of the equation (\ref{rho-e})
and integrating the resulting equation over $Q_t$, we have
\begin{align}
&\int_0^1\f{|\r^-|^2}{2}dx+\int_0^t\int_0^1
\big(\v+(\r^0\t^0)_\nu\big)|\r^-_x|^2dxd\tau
+\int_0^t\int_0^1\big(s\chi^\v(\sqrt{\t})|\r^-|^2
+s\chi^\v(p_s(\t))\r^-\big)dxd\tau
\nn\\
&~~+\int_0^t\big(\alpha^0|\r^-(0,\tau)|^2
+\alpha^0s\bar\r^0\r^-(0,\tau)\big)d\tau
+\int_0^t\big(\alpha^1|\r^-(1,\tau)|^2
+\alpha^1s\bar\r^1\r^-(1,\tau)\big)d\tau
\nn\\
&=-\int_0^t\int_0^1\r^-\r^-_x(\r^0_\v\t^0_x)_\v dxd\tau
\nn \\
&\leq\int_0^t\int_0^1\biggl(\frac{\|(\r^0_\v\t^0_x)_\v
\|_{L^\infty(Q_T)}}{2\v} |\r^-|^2+\frac{\v}{2}|\r^-_x|^2\biggl)
dxd\tau. \nn
\end{align}
Notice that $\r^-\geq 0$. Thus we have that
$$
\int_0^1|\r^-|^2dx \leq\frac{\|(\r^0_\v\t^0_x)_\v
\|_{L^\infty(Q_T)}}{2\v} \int_0^t\int_0^1|\r^-|^2 dxd\tau.
$$

 By Gronwall's inequality, we can see that
$\r^-\equiv 0$. Thus $\r=\r^+\geq0$. This and (3.8) imply that
$\r\in Y\hookrightarrow\hookrightarrow X$.

Similarly, by the $L^2$-theory of quasi-linear parabolic equations
\cite{LSU}, there exists a solution $\t\in W^{2,1}_2(Q_T)$ for the
system (\ref{theta-e})-(\ref{theta-b}) and
\begin{align}\label{rhomapping1}
\|\t\|_{W^{2,1}_2(Q_T)}\leq
C(\v^{-1},\nu^{-1},\|\r^0\|_X,\|\t^0\|_X,T).
\end{align}
Let $\t^+=\max\{\t,0\}$, $\t^-=\max\{-\t,0\}$. Then $\t=\t^+-\t^-$.
Multiplying $\t^-/(\r+\sigma)$ on both sides of the equation
(\ref{theta-e}) and integrating the resulting equation over $Q_t$,
we can get
$$
\begin{array}{l}
\di\int_0^1\f{|\t^-|^2}{2}dx+\int_0^t\int_0^1\f{\k^\v}{\r+\s}|\t_x^-|^2dxd\tau+\int_0^t\int_0^1\f{s(\l+\t)p_s(\t)}{(\r+\s)}\t^-dxd\tau\\
\di
+\int_0^t\int_0^1s\l\r\chi^\v(\sqrt{\t^0})\f{\t^-}{\r+\s}dxd\tau+\int_0^t\f{\t^-(1,\tau)}{\r(1,\tau)+\s}\b^1(s\bar\t^1+\t^-(1,\tau))d\tau\\
\di+\int_0^t\f{\t^-(0,\tau)}{\r(0,\tau)+\s}\b^1(s\bar\t^1+\t^-(0,\tau))d\tau=\int_0^t\int_0^1s\r\chi^\v(\sqrt{\t^0})\f{|\t^-|^2}{\r+\s}dxd\tau\\
\di+\int_0^t\int_0^1\k^\v\t_x^-\f{\r_x\t^-}{(\r+\s)^2}dxd\tau+\int_0^t\int_0^1\big[(\v+(\r^0\t^0)_\nu)\r_x
+\r(\r^0_\v\t^0_x)_\v\big]\t^-_x\f{\t^-}{\r+\s}dxd\tau.
\end{array}
$$
Since $p_s(\t)=0$ for $\t\leq0$, we observe that
$(\l+\t)p_s(\t)\t^-=0$ a.e in $\Omega_T$. By Cauchy inequality and
the estimations (3.8)-(3.9), we can estimate the terms in the
right hand side of the above equality. Thus we obtain
$$
\int_0^1|\t^-|^2dx \leq C(\v^{-1},\nu^{-1},\|\r^0\|_X,\|\t^0\|_X,T)
\int_0^t\int_0^1|\t^-|^2 dxd\tau.
$$
Gronwall's inequality gives that $\t^-\equiv 0$. Thus
$\t=\t^+\geq0$. This and (3.9) imply that $\t\in
Y\hookrightarrow\hookrightarrow X$.

We conclude that the mapping $M:X^2\times[0,1]\rightarrow X^2$ is a
compact mapping.

Now we prove the continuity of the mapping $M$. For any
$(\hat\r^0,\hat\t^0,\hat s)\in X^2\times[0,1]$, let
$(\hat\r,\hat\t)=M(\hat\r^0,\hat\t^0,\hat s)$. Then
\begin{equation}
 \hat\r_t-[(\v+(\hat\r^0\hat\t^0)_\nu)\hat\r_x+\hat\r(\hat\r^0_\v\hat\t^0_x)_\v]_x
+\hat s\hat\r\chi^\v(\sqrt{\hat\t^0})=\hat s\chi^\v(p_s(\hat\t^0)),
\label{rho-e1}
\end{equation}
\begin{equation}
\begin{array}{r}
(\hat\r+\s)\hat\t_t-(\hat\k^\v\hat\t_x)_x-\big[(\v+(\hat\r^0\hat\t^0)_\nu)\hat\r_x
+\hat\r(\hat\r^0_\v\hat\t^0_x)_\v\big]\hat\t_x\\[6pt]
-\hat s\hat\r\chi^\v(\sqrt{\hat\t^0})\hat\t+\hat s(\lambda+\hat\t)
p_s(\hat\t)\di=\hat s\lambda\hat\r\chi^\v(\sqrt{\hat\t^0}),
\end{array}
\label{theta-e1}
\end{equation}
with the initial and boundary conditions
\begin{equation}
\left\{
\begin{array}{lr}
\di(\v+(\hat\r^0\hat\t^0)_\nu)\hat\r_x+\hat\r(\hat\r^0_\v\hat\t^0_x)_\v
=\a^1(\hat s\bar\r^1-\hat\r),&\mbox{\rm at}~~~ x=1,
\\[3mm]
\di(\v+(\hat\r^0\hat\t^0)_\nu)\hat\r_x+\hat\r(\hat\r^0_\v\hat\t^0_x)_\v
=\a^0(\hat\r-\hat s\bar\r^0),&\mbox{\rm at}~~~
 x=0,\\[3mm]
\r(x,0)=\hat s\r_{0\v}(x),&\mbox{\rm for}~~~x\in\Omega.
\end{array}
\right. \label{rho-b1}
\end{equation}
and
\begin{equation}
\left\{
\begin{array}{lr}
\di \hat\k^\v\hat\t_x=\b^1(\hat s\bar\t^1-\hat\t),&\mbox{\rm at}~~~
 x=1,\\[3mm]
\di \hat\k^\v\hat\t_x=\b^0(\hat\t-\hat s\bar\t^0),&\mbox{\rm at}~~~
 x=0,\\[3mm]
\t(x,0)=\hat s\t_{0\v}(x),&\mbox{for}~~~x\in\Omega \, .
\end{array}
\right.\label{theta-b1}
\end{equation}

Denote $\tilde\r=\r-\hat\r$ and $\tilde\t=\t-\hat\t$. Then
$\tilde\r$ satisfies the following equation,
\begin{equation}
\begin{array}{l}
\di \tilde\r_t-(F-\hat F)_x+s\tilde\r\chi^\v(\sqrt{\t^0})+(s-\hat
s)\hat\r\chi^\v(\sqrt{\t^0})+\hat s\hat
\r[\chi^\v(\sqrt{\t^0})-\chi^\v(\sqrt{\hat\t^0})]\\
\di =(s-\hat s)\chi^\v(p_s(\t^0))+\hat
s[\chi^\v(p_s(\t^0))-\chi^\v(p_s(\hat\t^0))],
\end{array}
\label{density}
\end{equation}
where
$$
F=(\v+(\r^0\t^0)_\nu)\r_x+\r(\r^0_\v\t^0_x)_\v,
$$
$$
\hat
F=(\v+(\hat\r^0\hat\t^0)_\nu)\hat\r_x+\hat\r(\hat\r^0_\v\hat\t^0_x)_\v.
$$
Multiplying the equation (\ref{density}) by $\tilde\r$ and
 integrating over $Q_t$ gives
$$
\int_0^1\tilde\r^2(x,t) dx+\int_0^t\int_0^1\tilde\r_x^2
dxd\tau\leq C\left[\int_0^t\int_0^1\tilde\r^2dxd\tau+(s-\hat
s)^2+\|\r^0-\hat\r^0\|_X^2+\|\t^0-\hat\t^0\|_X^2\right]
$$
with
$C=C(\v^{-1},\nu^{-1},\|\r_{0\v}\|_{L^2(\Omega)},\|\r^0\|_X,\|\t^0\|_X,\|\hat\r^0\|_X,\|\hat\t^0\|_X,T)$.

Thus Gronwall inequality implies that
$$
\|\tilde\r\|^2_X\leq C\left[(s-\hat
s)^2+\|\r^0-\hat\r^0\|_X^2+\|\t^0-\hat\t^0\|_X^2\right].
$$
Similarly, we can derive the equation for $\tilde\t$ and get
$$
\|\tilde\t\|^2_X\leq C\left[(s-\hat
s)^2+\|\r^0-\hat\r^0\|_X^2+\|\t^0-\hat\t^0\|_X^2\right].
$$
Thus, the mapping $M:X^2\times[0,1]\rightarrow X^2$ is continuous.
The proof of Lemma 3.2 is complete. \quad
\endproof \vskip0.1in

In addition, for $s=0$ we can see that $M(\r,\t,0)=0$ for any
$(\r,\t)\in X^2$.  Thus, by the Leray--Schauder fixed point
theorem, there exists a fixed point for the mapping
$M(\cdot,\cdot,1):X^2\rightarrow X^2$ if all the functions
$(\r,\t)\in X^2$ satisfying
\begin{equation}\label{fixpointp}
(\r,\t)=M(\r,\t,s)
\end{equation}
for some $s\in[0,1]$ are uniformly bounded in $X^2$. In fact, by
the proof of Lemma \ref{mpcompconti}, $M$ maps $(\r,\t,s)\in
X^2\times[0,1]$ into $Y^2$. Therefore, if $(\r,\t)$ is a fixed
point of $M(\cdot,\cdot,1)$, then $(\r,\t)\in W^{2,1}_2(Q_T)$.

\begin{theorem}\label{exapps}{\it
Under the assumptions of Theorem \ref{thm1}, the system
(\ref{asys})-(\ref{asys-bi}) has a (strong) solution $(\r,\t)\in
W^{2,1}_2(Q_T)$ which satisfies
\begin{eqnarray}
&\displaystyle\rho\geq\underline{\r}_{\v,T}
\quad\mbox{and}\quad\underline{\t}_T\leq\t\leq
\overline{\t}_{T}\quad\mbox{for}\;\;\; (x,t) \in Q_T.
\label{positive}\\[5pt]
&\displaystyle\|\r\|_{L^\infty(I;L^4(\Omega))},\;\;
\|\r_x\|_{L^2(Q_T)},\;\;\|\r\r_x\|_{L^2(Q_T)}\;
\leq C_{\v,T},\nn\\[5pt]
&\displaystyle\|\t\|_{L^\infty(Q_T)},\;\;
\|\r\|_{L^\i(I;L^1(\Omega))},\;\; \| \t_x \|_{L^2(Q_T)}, \;\;\|
\rho_\v\t_x \|_{L^2(Q_T)}\;\leq C_{T} \, . \label{(3.17-2)}
\end{eqnarray}
where $\underline{\r}_{\v,T}$ and $C_{\v,T}$ are positive
constants which depend on $\v$ and $T$, independent of $\nu$;
$\underline{\t}_T$, and $\overline{\t}_{T}$ and $C_{T}$ are
positive constants, dependent upon $T$ and independent of $\v$ and
$\nu$. }
\end{theorem}

By the Leray-Schauder fixed point theorem, it suffices to prove
the uniform boundedness of functions $(\r,\t)\in X^2$ satisfying
the equation (\ref{fixpointp}) and (\ref{positive}).

\subsection{Uniform estimates}
We assume that $(\r,\t)\in X^2$ and therefore, $(\r,\t)=M(\r,\t,s)
\in Y^2$, for $s\in[0,1]$, $i.e.$, $(\r,\t)$ is a (strong)
solution of the following system,
\begin{align}
&\r_t-((\v+(\r\t)_\nu)\r_x)_x-(\r(\r_\v\t_x)_\v)_x
+s\r\chi^\v(\sqrt{\t})=s\chi^\v(p_s(\t)), \label{rho-ef}\\[5pt]
&(\r+\s)\t_t-(\k\t_x)_x-\left[(\v+(\r\t)_\nu)\r_x
+\r(\r_\v\t_x)_\v\right]\t_x\nn\\[2pt]
&~~~~~~~~~~-s\r\chi^\v(\sqrt{\t})\t+s(\l+\t)
p_s(\t)\di=s\lambda\r\chi^\v(\sqrt{\t}), \label{theta-ef}
\end{align}
with the initial and boundary conditions
\begin{equation}
\left\{
\begin{array}{lr}
\di(\v+(\r\t)_\nu)\r_x+\r(\r_\v\t_x)_\v
=\a^1(s\bar\r^1-\r),&\mbox{\rm at}~~~
x=1,\\[3mm]
\di(\v+(\r\t)_\nu)\r_x+\r(\r_\v\t_x)_\v
=\a^0(\r-s\bar\r^0),&\mbox{\rm at}~~~
 x=0,\\[3mm]
\r(x,0)=s\r_{0\v}(x),&\mbox{\rm for}~~~x\in\Omega,
\end{array}
\right. \label{rho-bf}
\end{equation}
and
\begin{equation}
\left\{
\begin{array}{lr}
\di \k^\v\t_x=\b^1(s\bar\t^1-\t),&\mbox{\rm at}~~~
 x=1,\\[3mm]
\di \k^\v\t_x=\b^0(\t-s\bar\t^0),&\mbox{\rm at}~~~
 x=0,\\[3mm]
\t(x,0)=s\t_{0\v}(x),&\mbox{for}~~~x\in\Omega,
\end{array}
\right.\label{theta-bf}
\end{equation}

In this subsection, we derive some uniform estimates for solutions
to the above initial-boundary value problems.

Firstly we add the equation (\ref{rho-ef}) multiplying by
$(\lambda+\t)$ into (\ref{theta-ef}) and then, integrate the
resulting equation over $Q_t$. We arrive at
\begin{equation}
\int_0^1(\l\r+\r\t+\s\t)(x,t)dx- \left. \int_0^tH_2(x,\tau) \right
|_{x=0}^{x=1}
d\tau\leq\int_0^1(\l\r_{0\v}+\r_{0\v}\t_{0\v}+\s\t_{0\v})(x)dx \nn
\end{equation}
where
\begin{equation}
H_2(x,\tau)=[\v\r_x+(\r\t)_\nu\r_x+\r(\r_\v\t_x)_\v](\l+\t)+\k^\v\t_x
\, . \nn
\end{equation}
With boundary conditions in (\ref{rho-bf})-(\ref{theta-bf}), we
have
\begin{eqnarray}
-H_2(x,\tau) \Big |_{x=0}^{x=1} &=&
\a^1(\r(1,\tau)-s\bar\r^1)(\l+\t(1,\tau))
+\a^0(\r(0,\tau)-s\bar\r^0)(\l+\t(0,\tau)) \nn \\ [2mm] && +
\b^1(\t(1,\tau)-s\bar\t^1)+\b^0(\t(0,\tau)-s\bar\t^0) \nn \\ [3mm]
&\geq & -\a^1s\bar\r^1\t(1,\tau)-\a^0s\bar\r^0\t(0,\tau) -\l
s(\a^1\bar\r^1+\a^0\bar\r^0)-s(\b^1\bar\t^1+\b^0\bar\t^0) \nn
\end{eqnarray}
and therefore,
\begin{equation}
\int_0^1(\l\r+\r\t+\s\t)(x,t)dx\leq
C_T+C\int_0^t\|\t(\cdot,\tau)\|_{C(\bar\Omega)}d\tau,
\label{(3.14)}
\end{equation}
where
$$
C_T=(\l+\|\t_{0\v}\|_{L^\i})\|\r_{0\v}\|_{L^1}+\s\|\t_{0\v}\|_{L^\i}
+\left[\l(\a^1\bar\r^1+\a^0\bar\r^0)+(\b^1\bar\t^1+\b^0\bar\t^0)\right]T.
$$

Similarly, subtracting the equation (\ref{theta-ef}) multiplied by
$\t^l/l$ from the equation (\ref{rho-ef}) multiplied by
$\t^{l+1}/(l+1)$ and integrating the resulting equation over $Q_t$,
we arrive at
\begin{eqnarray}\label{thetalp}
&& \int_0^1 (\r+\s)\t^{l+1}(x,t) dx- \left.
\int_0^tH_3(x,\tau)\right |_{x=0}^{x=1}d\tau +\int_0^t\int_0^1\k^\v
l(l+1)\t^{l-1}|\t_x|^2 dxd\tau
\nn\\
&& +s(l+1)\int_0^t\int_0^1(\lambda+\t)p_s(\t)\t^ldxd\tau =\int_0^1
(\r_{0\v}+\s)(\t_{0\v})^{l+1}(x) dx \label{rtk}
\\
&&+s\int_0^t\int_0^1\Big[l\t^{l+1} +\l(l+1)\t^l\Big]
\r\chi^\v(\sqrt{\t}) dxd\tau
+s\int_0^t\int_0^1\chi^\v(p_s(\t))\t^{l+1}dxd\tau \nn,
\end{eqnarray}
where
\begin{equation}
\begin{array}{lll}
\di -H_3(x,\tau) \Big |_{x=0}^{x=1}&=&\di
\a^1(\r(1,\tau)-s\bar\r^1)[\t(1,\tau)]^{l+1}
+\a^0(\r(0,\tau)-s\bar\r^0)[\t(0,\tau)]^{l+1}
\\[3mm]
&&\di+ (l+1)\b^1(\t(1,\tau)-s\bar\t^1)[\t(1,\tau)]^{l}+(l+1)
\b^0(\t(0,\tau)-s\bar\t^0)[\t(0,\tau)]^{l}\\[3mm]
&=&\di [\a^1\r(1,\tau)+(l+1)\b^1-\a^1s\bar\r^1]
[\t(1,\tau)]^{l+1}-(l+1)\b^1s\bar\t^1[\t(1,\tau)]^{l}\\[3mm]
&&\di+[\a^0\r(0,\tau)+(l+1)\b^0-\a^0s\bar\r^0]
[\t(0,\tau)]^{l+1}-(l+1)\b^0s\bar\t^0[\t(0,\tau)]^{l},
\\[3mm]
&\geq&\di -2^l(l+1)[\b^1(s\bar\t^1)^{l+1}+\b^0(s\bar\t^0)^{l+1}] \nn
\label{(3.43+)}
\end{array}
\end{equation}
when $l$ is large enough. Since $\t^l\leq \t^{1/2}+\t^{l+1}$ for
any $\t\geq0$ and $l\geq1$, by (\ref{(3.14)})-(\ref{thetalp}),
\begin{eqnarray}
\label{prertl} && \int_0^1 (\r+\s)\t^{l+1}(x,t) dx
+l(l+1)\int_0^t\int_0^1\k^\v \t^{l-1}|\t_x|^2
dxd\tau+sl\int_0^t\int_0^1(\lambda+\t)p_s(\t)\t^ldxd\tau
\nn\\
&&\leq C_{l,T}+C_0s\int_0^t\int_0^1l\big(1+\t^{l+1/2}\big) \r\t
dxd\tau
\nn\\
&& \leq C'_{l,T}+C_0sl\int_0^t \|\t(\cdot,\tau)
\|_{L^\infty(\Omega)}^{l+3/2}  d\tau
\end{eqnarray}
where
$$
C_{l,T}=\int_0^1 (\r_0+\s)\t_0^{l+1}(x) dx +
2^l(l+1)[\b^1(s\bar\t^1)^{l+1}+\b^0(s\bar\t^0)^{l+1}]
$$
and
$C_{l,T}'=C_{l,T}+Cl$. Recall the Gagliardo--Nirenberg inequality
$$
\|f\|_{L^\infty(\Omega)}\leq C
\|f\|_{L^2(\Omega)}+C\|f\|_{L^2(\Omega)}^{1/2}
\|f_x\|_{L^2(\Omega)}^{1/2},\quad \forall\,f\in H^1(\Omega).
$$
With $f=\t^{\frac{l+1}{2}}$ in the above inequality,
we obtain
\begin{equation}\label{ttttt}
\|\t(\cdot,\tau)\|_{L^\infty(\Omega)}^{l+3/2}\leq
\frac{C_2}{2}\int_0^1\t^{l+3/2}(x,\tau)dx
+C_1\big\|\t^{\f{l+1}{2}}(\cdot,\tau)\big\|_{L^2(\Omega)}^{\f{2l+3}{2l+2}}
\big\|(\t^{\f{l+1}{2}})_x(\cdot,\tau)\big\|_{L^2(\Omega)}^{\f{2l+3}{2l+2}}
\nn
\end{equation}
and by H\"{o}lder's inequality,
\begin{eqnarray}
&&\int_0^t\big\|\t^{\f{l+1}{2}}
(\cdot,\tau)\big\|_{L^2(\Omega)}^{\f{2l+3}{2l+2}}
\big\|(\t^{\f{l+1}{2}})_x(\cdot,\tau)
\big\|_{L^2(\Omega)}^{\f{2l+3}{2l+2}}d\tau\nn\\
&&\leq Cl\int_0^t\int_0^1(\t^\f{l+1}{2})^{\f{4l+6}{2l+1}}dxd\tau
+\frac{1}{(l+1)C_0C_1}\int_0^t\int_0^1\k|(\t^\f{l+1}{2})_x|^2dxd\tau
\nn\\
&&\leq Cl\int_0^t\int_0^1\t^\f{(l+1)(2l+3)}{2l+1}dxd\tau
+\frac{l+1}{4C_0C_1}\int_0^t\int_0^1\k\t^{l-1}|\t_x|^2dxd\tau. \nn
\end{eqnarray}
It follows that
\begin{eqnarray}
\int_0^t\|\t(\cdot,\tau)\|_{L^\infty(\Omega)}^{l+3/2}d\tau \leq
\frac{C_2}{2} \int_0^t\int_0^1\t^{l+3/2}dxd\tau
+Cl\int_0^t\int_0^1\t^\f{(l+1)(2l+3)}{2l+1}dxd\tau
\nn\\
+\frac{l+1}{4C_0}\int_0^t\int_0^1\k_1\t^{l-1}|\t_x|^2dxd\tau \nn
\, .
\end{eqnarray}
By the assumption (\ref{sat-a}), we observe that $C_0 C_2
\t^{l+\f32}\leq p_s(\t)\t^{l+1}+C$ for all $\t\geq0$. Substituting
the last inequality into (\ref{prertl}) gives
\begin{eqnarray}\label{prertf}
&&\int_0^1 (\r+\s)\t^{l+1}(x,t) dx
+\frac{l(l+1)}{2}\int_0^t\int_0^1\k^\v
\t^{l-1}|\t_x|^2dxd\tau+\frac{sl}{2}\int_0^t\int_0^1p_s(\t)\t^{l+1}dxd\tau
\nn\\
&&\leq C_{l,T}'+C_3sl^2
\int_0^t\int_0^1\t^\f{(l+1)(2l+3)}{2l+1}dxd\tau,
\end{eqnarray}
for $l$ being large enough. Let $l_0$ be a positive integer
satisfying
$$
\frac{(l_0+1)(2l_0+3)}{2l_0+1} =l_0+1 + \frac{2l_0+2}{2l_0+1}<l_0
+ 1 +(1+\eta)
$$
where $\eta$ is defined in (\ref{sat-a}). By noting the fact
$$
C_3 \t^\f{(l_0+1)(2l_0+3)}{2l_0+1} \le \frac{1}{4l_0} p_s(\theta)
\theta^{l_0+1} + (C l_0)^{l_1}
$$
with $l_1 = 2(l_0+2+\eta)/\eta$, we have
\begin{equation}\label{prertff}
\int_0^1 (\r+\s)\t^{l_0+1}(x,t) dx
+\frac{l_0(l_0+1)}{2}\int_0^t\int_0^1\k^\v \t^{l_0-1}|\t_x|^2dxd\tau
+\frac{sl}{4}\int_0^t\int_0^1p_s(\t)\t^{l_0+1}dxd\tau\leq
C_{l_0,T}'', \nn
\end{equation}
where $C_{l_0,T}''=C_{l_0,T}'+C_T (C l_0)^{l_1}$ for some constant
$C_T$ independent of $l_0$. Furthermore,
\begin{equation}\label{tff}
\sup_{0\leq t\leq
T}\int_0^1\t^{l_0+1}(x,t)dx+\int_0^T\int_0^1|(\t^{\frac{l_0+1}{2}})_x|^2
dxdt \leq C_{l_0,T}'' \nn
\end{equation}
and by the Sobolev embedding inequality,
\begin{equation}\label{tff2}
\int_0^T\|\t\|_{L^\infty(\Omega)}^{l_0+1} dxdt \leq
C_T^{l_0+1}C_{l_0,T}''. \nn
\end{equation}
Since $l_0$ is a fixed positive integer dependent solely upon
$\eta$, we obtain the estimate
\begin{equation}\label{tff3}
\int_0^T\|\t\|_{L^\infty(\Omega)} dxdt \leq C_T.
\end{equation}
From (\ref{(3.14)}) and (\ref{prertl}), we get
\begin{equation}\label{rrtl1}
\sup_{0\leq t\leq T}\int_0^1(\r+\r\t)dx\leq C_T
\end{equation}
and
\begin{eqnarray}\label{prertl22}
&& \int_0^1 (\r+\s)\t^{l+1}(x,t) dx +l(l+1)\int_0^t\int_0^1\k^\v
\t^{l-1}|\t_x|^2
dxd\tau+sl\int_0^t\int_0^1(\lambda+\t)p_s(\t)\t^ldxd\tau
\nn\\
&&\leq C_{l,T}+Csl\int_0^t\int_0^1\r\t
dxd\tau+Csl\int_0^t\|\t\|_{L^\infty(\Omega)}^{1/2}\int_0^1
\r\t^{l+1} dxd\tau
\nn\\
&&\leq (C_{l,T}+C_Tl)+Csl\int_0^t\|\t\|_{L^\infty(\Omega)}^{1/2}
\int_0^1(\r+\sigma)\t^{l+1} dxd\tau. \nn
\end{eqnarray}
Moreover, by using Gronwall's inequality,
$$
\int_0^1 (\r+\s)\t^{l+1}(x,t) dx\leq (C_{l,T}+C_Tl)+(C_{l,T}+C_Tl)
e^{C_T l}
$$
and
$$
\| \t \|_{L^{l+1}(Q_T)} \le [2(C_{l,T}+C_Tl)]^{\frac{1}{l+1}}
e^{C_T} \, .
$$
On the other hand, by taking $l\rightarrow\infty$, we have
\begin{equation}
\label{tlinfn} \|\t\|_{L^\infty(Q_T)}\leq C_T
\end{equation}
where we have noted the fact
$$
C_{l,T}^{\frac{1}{l+1}} \leq C_T \, .
$$
Moreover, by taking $l=1$ in the equation (\ref{prertl}), we
obtain
\begin{equation}
\di \int_0^1 (\r+\s)\t^{2}(x,t) dx+\f{1}{2}\int_0^t\int_0^1
\big((\k_1+\k_2|\r_\v|^2)|\t_x|^2+s\t^2 p_s(\t)\big) dxd\tau\leq C_T
\nn,
\end{equation}
which implies that
\begin{equation}
\| \t_x \|_{L^2(Q_T)},\;\;\| \rho_\v\t_x \|_{L^2(Q_T)}\leq C_{T}
\, \, . \label{(3.17)}
\end{equation}

Secondly we present some estimates for $\rho$. By multiplying $\r$
on both sides of the equation (\ref{rho-ef}) and integrating the
resulting equation over $Q_T$, with Gronwall's inequality we get
\begin{align}\label{rho2pv}
\sup_{0\leq t\leq T}\int_0^1\r^2dx+\int_0^T\int_0^1|\r_x|^2dxdt
\leq C_{\v,T}+C(\v,\|(\r_\v\t_x)_\v\|_{L^\infty(Q_T)}) \le
C_{\v,T},
\end{align}
which together with the Sobolev embedding inequality gives
$$
\int_0^T\int_0^1\r^6dxdt\leq C_{\v,T}.
$$

Once again, multiplying $\r^3$ on both sides of the equation
(\ref{rho-ef}) and integrating the resulting equation over $Q_T$
lead to
\begin{align}\label{rho2pv2}
\sup_{0\leq t\leq
T}\int_0^1\r^4dx+\int_0^T\int_0^1\r^2|\r_x|^2dxdt \leq C_{\v,T}.
\end{align}

From (\ref{tlinfn}), (\ref{(3.17)}) and (\ref{rho2pv}), we
conclude that $(\r,\t)$ is uniformly bounded in $X^2$. Thus, by
the Leray--Schauder fixed point theorem, there exists a fixed
point $(\r^{\v,\nu},\t^{\v,\nu})$ for the mapping
$M(\cdot,\cdot,1):X^2\rightarrow X^2$ and
$(\r^{\v,\nu},\t^{\v,\nu})$ is a solution of the system
(\ref{asys})-(\ref{asys-bi}).

\subsection{Positivity of the approximate solutions}\label{uniest}
Finally we prove the positivity of the approximate solutions
$(\r^{\v,\nu},\t^{\v,\nu})$. Let $\tilde\t^\delta=\t e^t-\delta$.
Then $\tilde\t^\delta$ is the solution of the following problem,
\begin{align}
&(\r+\s)\tilde\t_t^\delta-(\k^\v\tilde\t_x^\delta)_x
-\left[(\v+(\r\t)_\nu)\r_x +\r(\r_\v\t_x)_\v\right]\tilde\t_x^\delta
-(\r+\sigma)\tilde\t^\delta-\r\chi^\v(\sqrt{\t})
\tilde\t^\delta\nn\\[6pt]
&~+\tilde q(\t e^t,\delta)\tilde\t^\delta\di =\r\chi^\v(\sqrt{\t})\t
e^t+\l\r\chi^\v(\sqrt{\t})e^t+(\r+\sigma)\delta -(\l+e^{-t}\delta)
p_s(e^{-t}\delta)e^t, \label{theta-e3}
\end{align}
with the initial and boundary conditions
\begin{equation}
\left\{
\begin{array}{lr}
\di
\k^\v\tilde\t^\delta_x+\b^1\tilde\t^\delta=\b^1(\bar\t^1e^t-\delta),&\mbox{\rm
at}~~~ x=1,\\[3mm]
\di -\k^\v\tilde\t^\delta_x+\b^0\tilde\t^\delta
=\b^1(\bar\t^0e^t-\delta),&\mbox{\rm at}~~~ x=0,\\[3mm]
\tilde\t^\delta(x,0)=\t_{0\v}(x)-\delta,&\mbox{for}~~~x\in\Omega,
\end{array}
\right.\label{theta-b3}
\end{equation}
where
$$
\tilde q(\tilde\t,\delta)=\frac{(\l+e^{-t}\tilde\t)
p_s(e^{-t}\tilde\t)-(\l+e^{-t}\delta)
p_s(e^{-t}\delta)}{\tilde\t-\delta}e^t\geq0.
$$
By the assumption (\ref{sat-a}), the right hand side of the
equations (\ref{theta-e3})-(\ref{theta-b3}) are nonnegative if
$\delta$ is small enough (independent of $\v$ and $\nu$).
Multiplying $(\tilde\t^\delta)^-/(\r+\sigma)$ on both sides of the
equation (\ref{theta-e3}) and integrating the resulting equation
over $Q_t$, we derive $\tilde\t^\delta\geq 0$, i.e. $\t\geq
e^{-T}\delta$, which together with (\ref{tlinfn}) implies that
\begin{equation}\label{unibdtheta}
\underline{\t}_{T}\leq\theta(x,t)\leq\overline{\t}_{T}
\;\;\;\mbox{for}\;\;\; (x,t) \in Q_T.
\end{equation}
where $\underline{\t}_T$ and $\overline{\t}_{T}$ are positive
constants independent of $\v$ and $\nu$.

For $\r$, we define $\r^\delta=\r-\delta$. Then $\r^\delta$ is the
solution of the following equation
\begin{equation}\label{rho-e2}
 \r^\delta_t-((\v+(\r\t)_\nu)\r^\delta_x)_x-(\r^\delta(\r_\v\t_x)_\v)_x
+\r^\delta\chi^\v(\sqrt{\t})=\chi^\v(p_s(\t))+\delta[(\r_\v\t_x)_\v]_x
-\delta\chi^\v(\sqrt{\t}),
\end{equation}
with the initial and boundary conditions
\begin{equation}
\left\{
\begin{array}{lr}
\di(\v+(\r\t)_\nu)\r^\delta_x+\r^\delta(\r_\v\t_x)_\v+\a^1\r^\delta
=\a^1(\bar\r^1-\delta)-\delta(\r_\v\t_x)_\v,&\mbox{\rm
at}~~~ x=1,\\[3mm]
\di-(\v+(\r\t)_\nu)\r^\delta_x-\r^\delta(\r_\v\t_x)_\v
+\a^0\r^\delta=\a^0(\bar\r^0-\delta)+\delta(\r_\v\t_x)_\v,&\mbox{\rm
at}~~~ x=0,\\[3mm]
\r^\delta(x,0)=\r_{0\v}(x)-\delta,&\mbox{\rm for}~~~x\in\Omega.
\end{array}
\right. \label{rho-b2}
\end{equation}
Since $\chi^\v(p_s(\t))\geq\v$, the right hand side of the equations
(\ref{rho-e2})-(\ref{rho-b2}) are nonnegative if
$$
\delta=\min\biggl\{\frac{\v}{2},\;\;
\frac{\v}{1+2\|(\r_\v\t_x)_\v\|_{C^1(\overline Q_T)}}\biggl\},
$$
in which case $\r^\delta\geq0$, or eqivalently $\r\geq\delta$. On
the other hand, from (\ref{(3.17)}) we have
$$
\|(\r_\v\t_x)_\v\|_{C^1(\overline
Q_T)}\leq\frac{1}{\v^2}\|\r_\v\t_x\|_{L^1(Q_T)}\leq C_\v.
$$
Thus, there exists a positive constant $\underline{\r}_{\v,T}$
such that
\begin{equation}\label{lowerbdrho}
\rho\geq\underline{\r}_{\v,T}\;\;\;\mbox{for}\;\;\; (x,t) \in Q_T
\,.
\end{equation}

\section{Global existence} \setcounter{equation}{0}
We have constructed an approximate solution $(\rho^{\v, \nu},
\theta^{\v, \nu})$ to the system (\ref{asys2}) and (\ref{asys-bi})
(or equvilently (\ref{asys})-(\ref{asys-bi})) in the last section.
In this section, we prove the global existence of weak solutions
for the system (\ref{sys})-(\ref{sys-i}). Firstly we fix $\v>0$
and study the convergence as $\nu\rightarrow0$.

Since the system (\ref{rho-ef})-(\ref{theta-ef}) reduces to
(\ref{asys})-(\ref{asys-bi}) when $s=1$, the uniform estimates
(\ref{tlinfn}), (\ref{(3.17)}), (\ref{rho2pv}) and
(\ref{unibdtheta}) given in the last section still hold for the
approximate solution $(\rho^{\epsilon, \nu}, \theta^{\epsilon,
\nu})$. We rewrite the first equation in (\ref{asys}) by
$$
\r_t=-f_x+g
$$
with $g$ uniformly bounded in $L^2(Q_T)$ and
$$
f=(\v+(\r\t)_\nu)\r_x+\r(\r_\v\t_x)_\v.
$$
Since $\r$ is uniformly bounded in $L^\infty(I;L^2(\Omega))\cap
L^2(I;H^1(\Omega))\hookrightarrow L^6(Q_T)$, we derive that
$$\|\r(\r_\v\t_x)_\v\|_{L^2(Q_T)},\;\;\|(\r\t)_\nu\|_{L^6(Q_T)},\;\;
\|(\r\t)_\nu\r_x\|_{L^\frac{5}{4}(Q_T)}\; \leq C_{\v,T}
$$
and
$$
\|\r_t\|_{L^{5/4}(I;W_0^{-1,{5/4}}(\Omega))}\leq C_{\v,T}.
$$
From the first equation in (\ref{asys2}) we derive that
$$
\|(\r\t+\sigma\t)_t\|_{L^{5/4}(I;W_0^{-1,{5/4}}(\Omega))}\leq
C_{\v,T}
$$
where we have noted (\ref{tlinfn}), and moreover, from
(\ref{rho2pv}), we observe that $\r^{\v,\nu}$ is uniformly bounded
in $L^6(I;L^6(\Omega))\cap L^2(I;H^1(\Omega))$ and $\r^{\v,\nu}_t$
is uniformly bounded in $L^{5/4}(I;W_0^{-1,5/4}(\Omega))$. Using
Aubin--Lions lemma, we conclude that there exists a sequence
$\nu_j\rightarrow0$ such that
\begin{equation}
\begin{array}{ll}
&\r^{\v,\nu_j}\rightarrow \r^\v~~ \mbox{\rm strongly~~
in}~~~L^p(Q_T)~(\forall 1\leq p<6), \\ [2mm]
&\r^{\v,\nu_j}\rightarrow \r^\v~~ \mbox{\rm strongly~~
in}~~~L^2(I;C(\overline\Omega)),\\ [2mm]
&\r^{\v,\nu_j}{\rightharpoonup} \r^\v ~~~ \mbox{\rm weakly ~~
in}~~~L^2(I,H^1(\Omega)),\\ [2mm]
&\r^{\v,\nu_j}_t{\rightharpoonup} \r^\v_t ~~~ \mbox{\rm weakly ~~
in}~~~L^{5/4}(I;W_0^{-1,{5/4}}(\Omega))
\end{array}
\label{convnurho}
\end{equation}
and
$$
\r^{\v,\nu_j}(0,\cdot)\rightarrow
\r^\v(0,\cdot)\quad\mbox{and}\quad\r^{\v,\nu_j}(1,\cdot)\rightarrow
\r^\v(1,\cdot)\quad\mbox{strongly \;in}\;\;\;L^2(0,T).
$$
Similarly, by noting the uniform estimates (\ref{tlinfn}),
(\ref{(3.17)}) and (\ref{unibdtheta}), we conclude that there
exists a subsequence of $\t^{\v,\nu_j}$ (also denoted by
$\t^{\v,\nu_j}$) such that
\begin{equation}
\begin{array}{ll}
&\t^{\v,\nu_j}\rightarrow \r^\v~~ \mbox{\rm strongly~~
in}~~~L^p(Q_T)~(\forall 1\leq p<\infty), \\ [2mm]
&\t^{\v,\nu_j}\rightarrow \t^\v~~ \mbox{\rm strongly~~
in}~~~L^2(I;C(\overline\Omega)),\\ [2mm]
&\t^{\v,\nu_j}{\rightharpoonup} \t^\v ~~~ \mbox{\rm weakly ~~
in}~~~L^2(I,H^1(\Omega)),\\ [2mm]
&(\r^{\v,\nu_j}\t^{\v,\nu_j}+\sigma\t^{\v,\nu_j})_t{\rightharpoonup}
(\r^\v\t^\v+\sigma\t^\v)_t ~~~ \mbox{\rm weakly ~~
in}~~~L^{5/4}(I;W_0^{-1,{5/4}}(\Omega))
\end{array}
\label{convnutheta}
\end{equation}
and
$$
\t^{\v,\nu_j}(0,\cdot)\rightarrow
\t^\v(0,\cdot)\quad\mbox{and}\quad\t^{\v,\nu_j}(1,\cdot)\rightarrow
\t^\v(1,\cdot)\quad\mbox{strongly \;in}\;\;\;L^p(0,T),\;\;1\leq
p<\infty.
$$
Since $(\r^{\v,\nu_j},\t^{\v,\nu_j})$ is a strong solution of the
system (\ref{asys2}) and (\ref{asys-bi}), it satisfies
\begin{align*}
&\int_0^T\alpha^0(\r^{\v,\nu_j}(0,t)-\bar\r^0)\phi(0,t)dt
+\int_0^T\alpha^1(\r^{\v,\nu_j}(1,t)-\bar\r^1)\phi(1,t)dt
\nn\\
&~~+\int_0^T\int_\Omega\rho_t^{\v,\nu_j}\phi dxdt
+\int_0^T\big[(\v+(\r^{\v,\nu_j}\t^{\v,\nu_j})_\nu)\r^{\v,\nu_j}_x
+\r^{\v,\nu_j}(\r^{\v,\nu_j}_\v\t^{\v,\nu_j}_x)_\v\big]\phi_xdxdt
\\
&=\int_0^T\int_\Omega\chi^\v(p_s(\t^{\v,\nu_j}))dxdt
-\int_0^T\int_\Omega\r^{\v,\nu_j}\chi^\v(\sqrt{\t^{\v,\nu_j}})\phi
dxdt
\end{align*}
and
\begin{align*}
&\int_0^T\int_0^1[(\r^{\v,\nu_j}+\s)\t^{\v,\nu_j}]_t\psi dxdt
+\int_0^T\beta^0(\t^{\v,\nu_j}(0,t)-\bar\t^0)\psi(0,t)dt
+\int_0^T\beta^1(\t^{\v,\nu_j}(1,t)-\bar\t^1)\psi(1,t)dt
\\[6pt]
&~~+\int_0^T\alpha^0(\r^{\v,\nu_j}(0,t)-\bar\r^0)\t^{\v,\nu_j}(0,t)\psi(0,t)dt
+\int_0^T\alpha^1(\r^{\v,\nu_j}(1,t)-\bar\r^1)\t^{\v,\nu_j}(1,t)\psi(1,t)dt
\nn\\
&~~+\int_0^T\int_0^1\k^\v\t^{\v,\nu_j}_x\psi_x dxdt
+\int_0^T\int_0^1\big[(\v+(\r^{\v,\nu_j}\t^{\v,\nu_j})_\nu)
\r^{\v,\nu_j}_x\t^{\v,\nu_j}+\r^{\v,\nu_j}(\r^{\v,\nu_j}_\v
\t^{\v,\nu_j}_x)_\v\t^{\v,\nu_j}\big]\psi_x dxdt
\\[6pt]
&~~+\int_0^T\int_0^1(\l+\t^{\v,\nu_j}) p_s(\t^{\v,\nu_j})\psi\big]
dxdt
\\[6pt]
&=\lambda\int_0^T\int_0^1\r^{\v,\nu_j}\chi^\v(\sqrt{\t^{\v,\nu_j}})\psi
dxdt
+\lambda\int_0^T\int_0^1\t^{\v,\nu_j}\chi^\v(p_s(\t^{\v,\nu_j}))\psi
dxdt,
\end{align*}
for any $\phi,\psi \in L^5(I;W^{1,5}(\Omega))$. By taking the
limit $j\rightarrow\infty$, we obtain a global weak solution
$(\r^\v,\t^\v)$ to the approximate system
\begin{align}
&\r_t-((\v+\r\t)\r_x)_x-(\r(\r_\v\t_x)_\v)_x=-\G_\v, \nn
\\[3mm]
&(\r\t+\s\t)_{t}-(\kappa^\v\t_x)_x-((\v+\r\t))
\r_x\t)_x-(\r(\r_\v\t_x)_\v\t)_x \label{asys3}
\\
&=\l\G_\v+(\l+\t) \left ( \chi^\v(p_s(\t))-p_s(\t)\right ) , \nn
\end{align}
with the boundary and initial conditions
\begin{equation}
\begin{array}{l}
\di (\v+\r\t)\r_x+\r(\r\t_x)_\v\big|_{x=1}=\a^1(\bar\r^1-\r(1,t)),\\[3mm]
\di (\v+\r\t)\r_x+\r(\r\t_x)_\v\big|_{x=0}=\a^0(\r(0,t)-\bar\r^0),\\[3mm]
\di \r(x,0)=\r_{0\v}(x):=\r_0\ast\eta_\v(x)+\v,\\[3mm]
\di \k^\v\t_x|_{x=1}=\b^1(\bar\t^1-\t(1,t)),\\[3mm]
\di \k^\v\t_x|_{x=0}=\b^0(\t(0,t)-\bar\t^0),\\[3mm]
\di \t(x,0)=\t_{0\v}(x):=\t_0\ast\eta_\v(x).
\end{array}
\label{asys-bi3}
\end{equation}

Secondly, we study the convergence as $\v\rightarrow0$.  To take
the limit $\,\v\rightarrow0$, we need more uniform estimates for
$\rho$ with respect to $\v$.

Clearly the system (\ref{asys})-(\ref{asys-bi}) reduces to the
system (\ref{asys3})-(\ref{asys-bi3}) when $\nu=0$. Then the
uniform estimates (\ref{positive}) and (\ref{(3.17-2)}) hold for
the obtained solution $(\r^\v,\t^\v)$. From (\ref{rho2pv2}) we see
that
$$
\|\r\t\r_x\|_{L^2(Q_T)}\leq C_{\v,T}
$$
and from the first equation of (\ref{asys3}) we deduce that $\r_t\in
L^2(I;H^{-1}_0(\Omega))$. Note that $\ln\r\in L^2(I;H^1(\Omega))$.
By multiplying the first equation of (\ref{asys3}) by $\ln\r$ and
integrating the equation over $Q_t$, we arrive at
$$
\begin{array}{l}
\di \int_0^1 \r\ln\r(x,t)dx-\int_0^1 \r(x,t)
dx+\int_0^t[\v\r_x+\r\t\r_x +\r(\r_\v\t_x)_\v]\ln\r\Big
|_{x=0}^{x=1}d\tau +\int_0^t\int_0^1\t\r_x^2dxd\tau
\\ [4mm]
\di\leq \int_0^1 \r_{0\v}\ln\r_{0\v}(x)dx-\int_0^1
\r_{0\v}dx-\int_0^t\int_0^1(\r_\v\t_x)_\v\r_xdxd\tau
-\int_0^t\int_0^1 (\r\sqrt{\t}- p_s(\t))\ln\r\, dxd\tau.
\end{array}
$$
Since
\begin{equation}
\begin{array}{ll}
\di \int_0^t\int_0^1|(\r_\v\t_x)_\v\r_x|dxd\tau&\di \leq
\f12\int_0^t\int_0^1\t\r_x^2dxd\tau
+\f12\int_0^t\int_0^1\f{|(\r_\v\t_x)_\v|^2}{\t}dxd\tau
\\ [4mm]
&\di \leq
\f12\int_0^t\int_0^1\t\r_x^2dxd\tau+C_T\|\r_\v\t_x\|_{L^2(Q_T)}^2
\\ [4mm]
&\di \leq \f12\int_0^t\int_0^1\t\r_x^2dxd\tau+C_T \, ,
\\ [4mm]
\end{array}
\nn \label{(3.56)}
\end{equation}
we get
\begin{equation}
\begin{array}{l}
\di \int_{[0,1]\cap\{\r\geq 1\}}
\r\ln\r(x,t)dx+\f12\int_0^t\int_0^1\t\r_x^2dxd\tau
+\int\int_{[0,1]\times[0,t]\cap\{\r\geq 1\}}\r\ln\r dxd\tau
\\ [4mm]
\di  \leq\int_0^1 \r_{0\v}|\ln\r_{0\v}|(x)dx
+\int_{[0,1]\cap\{\r\leq 1\}} \r|\ln\r|(x,t)dx
\\ [4mm]
\di\quad +\int\int_{[0,1]\times[0,t]\cap\{\r\leq 1\}}\r|\ln\r|
dxd\tau+\int\int_{[0,1]\times[0,t]\cap\{\r\geq 1\}} p_s(\t) \ln\r
dxd\tau+C_T
\\ [6mm]
\leq \di C_T,
\end{array}
\nn \label{(3.55)}
\end{equation}
which, together with (\ref{unibdtheta}), leads to
\begin{equation}
\|\r\ln\r\|_{L^\i(0,T;L^1(\Omega))},\quad \|\r_x\|_{L^2(Q_T)}\leq
C_T \, . \label{(3.24)}
\end{equation}
From the inequalities (\ref{rrtl1}) and (\ref{(3.24)}) we derive
\begin{equation}
\|\r\|_{L^2(0,T;H^1(\Omega))}\leq C_T\label{(3.25)}
\end{equation}
and
\begin{align}
\|\r\|_{L^\infty(\Omega)}^3 &\leq\|\r\|_{L^1(\Omega)}^3
+\|\r\|_{L^2(\Omega)}^{3/2}\|\r_x\|_{L^2(\Omega)}^{3/2}
\nn\\
&\leq
C_T+C\|\r\|_{L^1(\Omega)}^{3/4}\|\r\|_{L^\infty(\Omega)}^{3/4}\|\r_x\|_{L^2(\Omega)}^{3/2}
\nn\\
&\leq
C_T+\frac{1}{2}\|\r\|_{L^\infty(\Omega)}^3+C_T\|\r_x\|_{L^2(\Omega)}^2,
\nn
\end{align}
which results in
$$
\int_0^T\|\r\|_{L^\infty(\Omega)}^3dt\leq
C_T+C_T\int_0^T\|\r_x\|_{L^2(\Omega)}^2dt\leq C_T.
$$
Moreover, we have
\begin{equation}\label{ur4}
\int_0^T\int_0^1\r^4dxdt\leq
\biggl(\int_0^T\|\r\|_{L^\infty(\Omega)}^3dt\biggl)\biggl(\sup_{0\leq
t\leq T}\int_0^1\r dx\biggl)\leq C_T,
\end{equation}
i.e. $\r$ is uniformly bounded in $L^4(Q_T)$.

Finally, we let
\begin{equation}
B_1=H^1(\Omega),~B_2=L^4(\Omega),~B_3=W_0^{-1,6/5}(\Omega). \nn
\label{(3.32)}
\end{equation}
Then $B_1\hookrightarrow\hookrightarrow B_2\hookrightarrow B_3$
and $\{\r^\v\}$ is uniformly bounded in $L^4(I;B_2)\cap
L^2(I;B_1)$. From the first equation in (\ref{asys}), i.e.
$$
\r_t=\Big[\v\r_{x}+\r\t\r_x +\r(\r\t_x)_\v\Big]_x
-\r\chi^\v(\sqrt{\t})+\chi^\v(p_s(\t)) \,,
$$
we observe that $\{\r^\v_t\}$ is uniformly bounded in
$L^{6/5}(I;B_3)$. By Aubin--Lions lemma, $\{\r^\v\}$ is relatively
compact in $L^p(I;L^4(\Omega))$ for $(1\leq p<4)$. Thus, there
exists a sequence $\r^{\v_j}$ such that
$\lim_{j\rightarrow\infty}\v_j=0$ and
\begin{equation}
\begin{array}{ll}
&\r^{\v_j}\rightarrow \r~~ \mbox{\rm strongly~~
in}~~~L^p(I,L^4(\Omega))~(\forall 1\leq p<4), \\ [2mm]
&\r^{\v_j}\rightarrow \r~~ \mbox{\rm strongly~~
in}~~~L^2(I,C(\overline\Omega)), \\ [2mm]
&\r^{\v_j}{\rightharpoonup} \r ~~~ \mbox{\rm weakly ~~
in}~~~L^2(I,H^1(\Omega)),\\[2mm]
&\r_t^{\v_j}{\rightharpoonup} \r_t ~~ \mbox{\rm weakly ~~
in}~~~L^{6/5}(I;W_0^{-1,6/5}(\Omega)).
\end{array}
\label{(3.27)}
\end{equation}
Similarly, by (\ref{positive}) and (\ref{(3.17-2)}), there exists
a subsequence of $\t^{\v_j}$ (also denoted by $\t^{\v_j}$) such
that
\begin{equation}
\begin{array}{ll}
&\t^{\v_j}\rightarrow \t~~ \mbox{\rm strongly~~
in}~~~L^p(Q_T)~(\forall 1\leq p<\i), \\ [2mm]
&\t^{\v_j}\rightarrow\t~~ \mbox{\rm strongly~~
in}~~~L^2(I,C(\overline\Omega)), \\ [2mm]
&\t^{\v_j}{\rightharpoonup} \t ~~~ \mbox{\rm weakly ~~
in}~~~L^2(I,H^1(\Omega)),\\[2mm]
&(\r^{\v_j}\t^{\v_j}+\sigma\t^{\v_j})_t{\rightharpoonup}
(\r\t+\sigma\t)_t ~~ \mbox{\rm weakly ~~
in}~~~L^{6/5}(I;W_0^{-1,6/5}(\Omega)).
\end{array}
\label{(3.28)}
\end{equation}
Now we take the limit $j\rightarrow\infty$ and by (\ref{(3.27)})
and (\ref{(3.28)}), we obtain a weak solution $(\r,\t)$ which
satisfies (\ref{rdefeq}) and (\ref{tdefeq}). \quad \vskip0.1in

{\bf Acknowledgements} The authors wish to thank Professors P.
Lei, T. Yang, G. Yuan and X. Xu for helpful discussions.

\end{document}